\documentclass{amsart}
\usepackage{yfonts}
\newtheorem{theorem}{Theorem}

\newtheorem{lemma}{Lemma}
\newcommand{\RNum}[1]{\uppercase\expandafter{\romannumeral #1\relax}}
\newtheorem{corollary}{Corollary}
\newcommand{\beq}{\begin{equation}}
\newcommand{\eeq}{\end{equation}}

\begin{document}

\begin{center}
{\Large On Some Identities of Barred Preferential Arrangements}\\ 
 \vspace{5mm}
{\large S.Nkonkobe},\,V.Murali,\, \\
\vspace{2mm}\

{\it \footnotesize Department of Mathematics (Pure \& Applied)\\ Rhodes
University \\Grahamstown 6140 South Africa\\  snkonkobe@yahoo.com, v.murali@ru.ac.za}    \\
\vspace{2mm}
\end{center}

\section{Abstract}
%\begin{abstract}

A preferential arrangement of a finite set is an ordered partition. Associated with each such ordered partition is a chain of subsets or blocks endowed with a linear order. The chain may be split into sections by the introduction of a vertical bar, leading to the notion of a barred preferential arrangements. In this paper we derive some combinatorial identities satisfied by the number of possible barred preferential arrangements of an $n$-element set. We illustrate with some suitable examples highlighting some important consequences of the identities.  
%\end{abstract}

Mathematics Subject Classifications:05A18,05A19,05A16, 2013
%\section{Abstract}
\section{introduction}

The study of preferential arrangements in Combinatorics goes back to Cayley in \cite{cayley:1859}. Further studies of these objects were made by Gross\cite{gross:1962}, Mendelson\cite{mendelson:1982} and others. More recently Pippenger studied special kind of preferential arrangements and called them barred preferential arrangements in \cite{Pippenger:2010}. The same further generalised in \cite{barred:2013}. In this paper we derive several identities using combinatorial arguments, thus answering a comments raised by Pippenger in \cite{Pippenger:2010}.In preliminaries  we gather known results on preferential and barred preferential arrangements. In section 4 we state and prove several identities involving the number of barred preferential arrangements of an $n$-element set. In proving the identities we use combinatorial arguments. In section 5 we discuss some identities involving restricted barred preferential arrangements. Throughout the paper we illustrate the concept with some simple examples.
\section{Preliminaries}
In this section we collect all the preliminaries and known results setting up convenient notations for the discussion of barred preferential arrangements. We refer to the papers \cite{barred:2013} and \cite{Pippenger:2010} for discussions and results on barred preferential arrangements.\\

\noindent { {$ \mathbf 1^ \circ $ } Preferential arrangement. }\\
A partition of a set $X_n = \{1,2, \ldots,n\}$ with $n$ elements is a collection of disjoint non-empty subsets whose union is $X_n$. The subsets forming the partition are called blocks. If the blocks are arranged in a linear order, the ordered partition is referred as a preferential arrangement of $X_n$. see \cite{gross:1962}, \cite{mendelson:1982}. For example,
$\{\{1\}, \{2,3\} \}$ is a preferential arrangement of $X_3$ with two blocks with the natural ordering\\ $\{1\} \leq \{2,3\} $. For convenience we write this preferential arrangement as $1 \,\,\, 23$, reading the block containing 1 comes first and the block containing 23 comes second. From now onwards we will use the natural ordering of blocks as illustrated here. A warning: The natural ordering of the numbers $ 1,2, \ldots, n$ in not relevant here but the ordering of the blocks are. Let $Q_n$ and $J_n$ denote the set of preferential arrangements of $X_n$ and the number of preferential arrangements of $X_n$, so that $|Q_n|= J_n$. The numbers $J_n$ for $n \in \mathbb N_0$ are also known as Fubini numbers and a list of values of $J_n$ for $n= 0,1,2, \ldots$ can be found in OEIS sequence number A000670, \cite{oeis}. The number $J_n$ is interpreted as the number of outcomes in a race with $n$ participants assuming all of them finish the race and ties are allowed \cite{mendelson:1982} or the number of restricted preferential fuzzy subsets of $X_n$, \cite{murali:Nelsongenerating function}. The numbers $J_n$ satisfy a recurrence relation with the initial condition $J_0=1$, or identity as a finite sum, see \cite{mendelson:1982} or an infinite sum, see \cite{gross:1962} or a closed form involving Stirling numbers of the second kind \cite{Pippenger:2010}, respectively
\beq
J_{n+1} = \sum_{s=0}^n {n+1 \choose s} J_s, \quad or \quad J_n = \frac{1}{2} \sum_{s=0}^\infty s^n 2^{-s}, \quad or \quad J_n = \sum_{s=0}^n {n \brace s} s!\,.
\eeq
\noindent { {$ \mathbf 2^ \circ $ } Barred Preferential arrangements. }\\
The idea of introducing bars in between the blocks of a preferential arrangement seems first to appear in \cite{Pippenger:2010} even though the numbers associated with counting barred preferential arrangements appear much earlier in the literature \cite{com:74}. For example, for the preferential arrangement  $ \textgoth{X} = 1\,\,\,23$ considered in $ \mathbf 1^ \circ $, we can insert one bar $|$ in three different places in three different ways ( what a coincidence!) as follows: $1\, |\,\, 23 $ or $ |\,1 \,\,23 $ or $1 \,\, 23 \, |$, whereas a two-bar insertion into the same $\textgoth{X}$ would give rise to 6 barred preferential arrangements. viz., $1\, |\,|\,\, 23 $, $ |\,|\,1 \,\,23 $, $1 \,\, 23 \, |\,|$, $1\, |\,\, 23 \, |$, $ |\,1 \,\,23\, | $ or $|\,1 \,\,| 23 $. We now illustrate what we mean by sections associated with bars. Consider two 3-bar preferential arrangements of $X_7$ with 4 blocks.
\[\,\, \underbrace {7 \, \, 43 \,\, 1}_{1^{st}section} \,\,| \,\,\underbrace { \, \, 652 \,\, }_{2^{nd}section} \,\, | \,\,\underbrace { \,\, \,}_{3^{rd}section} | \,\,\underbrace { \,\, \,}_{4^{th}section}\,\, \]  

\[\,\, \underbrace { \, \, }_{1^{st} section} \,\,| \,\,\underbrace { \, \, 652 \,\, 1 }_{2^{nd}section} \,\, | \,\,\underbrace { \,\, 7 \,\, 43\,}_{3^{rd}section} | \,\,\underbrace { \,\, \,}_{4^{th}section}\,\, \]  
We notice that in general $m$ bars induce $m+1$ sections with some of the sections possibly empty.\\ 

In general we denote the number of barred preferential arrangements of $X_n$ with $m$ bars ($m \geq 0$) by $J_n^m$. In addition to the numbers $J_n^m$, we use $Q_n^m$ for the set of all barred preferential arrangements of $X_n$ with $m$ bars ($m \geq 0$). Hence $|Q_n^m| = J_n^m$. When $m=0$, $J_n^0$ ($|Q_n^0|$) identified as simply the number ( the set) of preferential arrangements without any bars and therefore is equal to $J_n$ ( $Q_n$ ) according to the notation in subsection {$ \mathbf 1^ \circ $ }. For instance $J_2^2= 15, J_4^5= 5340$ can be verified by direct manual counting of such barred preferential arrangements. The following closed form for $J_n^m$ and the recurrence relation with the initial condition $J_0^m = 1, J_1^m= m+1$ satisfied by it are obtained by Ahlbach et al in \cite{barred:2013}
\beq
J_n^m  = \sum_{s=0}^n {n \brace s} s! {m+s \choose m}, \quad J_n^m = \sum_{s=0}^n {n \choose s}J_s^0 \, J_{n-s}^{m-1} \,\,\, (m \geq 1)  
\eeq
The infinite sum representation for $J_n^m$ corresponding to similar expression for preferential arrangement is bit more involved. 
\beq
J_m^n = \frac{1}{2^{m+1} \, m!} \sum_{s=0}^\infty \frac{(s+1)^{\bar m} s^n} {2^s} \mbox{ where } (s+1)^{\bar m} = \sum_{t=1}^m  s(m,t) \times (s+1)^m
\eeq
again where $s(m,t)$ is the Stirling numbers of the first kind.\\

\noindent { {$ \mathbf 3^ \circ $ } Restricted Barred Preferential arrangements. }\\
We also discuss some identities based on barred preferential arrangements with some restrictions. Consider the following restriction on a preferential arrangement $\textgoth{X}$ of $X_n$ with $m$ bars. There are $m+1$ sections. We fix a specific section, say the ${(m+1)}^{th}$ section, namely the last section. After fixing the section, we impose the restriction that every other section to have at most one block. We then preferentially arrange the elements of the blocks of the fixed section. we then call such a barred preferential arrangement a restricted barred preferential arrangement. We denote the set of all restricted barred preferential arrangements of $X_n$ having $m$ bars by $H_n^m$. The order of the set $H_n^m$ is denoted by the number $I_n^m$. In the proofs involving the set $H_n^m$, it is occasionally useful to consider one  extra bar $\overset{*}{{|}}$ introduced on the elements of $H_n^m$ at the far right hand end of the last section. It will be clear in the proofs that such an introduction of an extra bar does not affect counting. Therefore we use set $*\P_n^m$ for elements of $H_n^m$ with  $\overset{*}{{|}}$. Both $H_n^m$ and $*\P_n^m$ have the same cardinality $I_n^m$.

% $ \ell \vdash n$ for the partition $ \ell = (l_1, l_2, \cdots, l_r) $ whose parts are the $ l_i$'s.

\section{The identities}
In this section we develop a number of interesting identities satisfied by $J_n^m$ and for various non-negative integral values of $n$ and $m$. The proofs are combinatorial arguments based on preferential arrangements. The first identity was found in \cite{Pippenger:2010} without proof or reference. Further the author of that paper states that such an identity can be proved combinatorially or otherwise. Here below we state and prove that result as a lemma.
\begin{lemma}\label{lemma:1}
\mbox{Let} $n$ \mbox{be an non-negative integer}. Then $J^0_{n+1}=\sum\limits_{s=0}^{n}\binom{n}{s}J^1_{s}$
\end{lemma}
\textbf{proof}:
What we are required to prove is that the number of preferential arrangements with no bars can be written as a sum of barred preferential arrangements with one bar.\\ 
The left hand side of the identity can be interpreted combinatorially as outcomes in a race of $n+1$ people ( assuming that all people finish the race )\cite{mendelson:1982}. The terms on the right hand side would correspond to races of $n$ people, of which the people are preferentially arranged with respect to a single bar ( hence the superscript of $J$ is 1 ). In proving the lemma, use a similar method to the one used in \cite{mendelson:1982} in establishing equation (7)). We need to associate each term on the left with one term on the right and visa versa for the identity to hold. Therefore we first mark a person out of $n+1$ people who are to run a race. We interpret the block having the marked person as a bar. Clearly all outcomes of the race would not be affected whether a specific person is marked or not. In our notation the number of outcomes is captured by $J_{n+1}^0$ (by definition of $J_{n+1}^0$). We ask the following question in order to identify each outcome that is counted by $J^0_{n+1}$is interpreted as an outcome in the sum on the right hand side and visa versa: How many people do not finish with the marked person?  The block having the marked person would serve as the bar in a barred preferential arrangement counted on the right hand side of the identity. There are $\binom{n}{s}$ ways of choosing $s$ people who do not finish with the marked person for each $s$ such that $0\leq s\leq n$. For each such $s$, there are $J^1_s$ ways of preferentially arranging the $s$ people with the block having the marked person serving as the bar. Therefore there are $\binom{n}{s} \times J^1_s$ many number of barred preferential arrangement of the $s$ people each of which is a preferential arrangement of $n+1$ persons finishing a race. Clearly the preceding argument uniquely identifies each outcome of the LHS as an outcome of the RHS.\\

%as a bar in a preferential arrangement of $s$ people, in which the preferential arrangement has one bar.

%Finding the total number of selecting the $s$ elements, preferentially arranging them in $J_{s}^{1}$ ways and summing over $s$ we obtain the right hand side of the identity.
% we need to define a bijective function.   

We generalize Lemma \ref{lemma:1} to the case of multiple bars in the following way, 
\begin{theorem}\label{theorem:1} For two integers $m,n\geq0$, $J^m_{n+1}\,=\,(m+1)\sum\limits_{s=0}^n\binom{n}{s} J^{m+1}_s$\quad  
\end{theorem}

\textbf{proof}: 
In forming a barred preferential arrangement of an $n+1$-persons as in Lemma one above, with $m$ bars, we mark one of the $n$ persons. The position of the marked person would be determined by two quantities, one: the section in which the marked person finishes in and two: the block in which he finishes within the section.  The block in which the marked person is in, is interpreted as the ${m+1}^{th}$ bar on the right hand side. Let us assume that the marked person is in section $i$. The question is: how many people do not finish with the marked person? Let us say there are $s$ of them, those $s$ people can be selected in $n\choose s$ ways. There are  $J^{m+1}_s $  ways of preferentially arranging the $s$ persons, where the ${m+1}^{th}$ bar is the block having the marked person. Summing over $s$ we have $\sum\limits^{n}_{s=0}{n\choose s} J^{m+1}_s$ number of preferential arrangements. If the marked person finishes in section $k$ instead of section $i$ above, we would still have  $\sum\limits^{n}_{s=0}{n\choose s} J^{m+1}_s$ many preferential arrangements.  Since there are $m+1$ sections hence the total number of preferential arrangements is $(m+1) \sum\limits^{n}_{s=0}{n\choose s} J^{m+1}_s$.

After putting the sums in  theorem \ref{theorem:1} in a Pascal like triangle for fixed values of $m$ and $n$, we obtain the following corollaries. 
\begin{corollary}
$J^1_2=4\times J^1_1$
\end{corollary}
We prove the corollary combinatorially by establishing a bijection between a set $Q^1_2$ whose cardinality is $J^1_2$ and a set $\{0,1\}_2\times\{0,1\}_1\times Q_1^1$ whose cardinality is $4\times J^1_1$. 
We denote by  $Q^m_n$ the set of barred preferential arrangements of a n-element set having $m$ bar. We have, $Q^1_2=\{|ab, ab|,  a\thinspace\thinspace\thinspace b|, b\thinspace\thinspace\thinspace a|, |a\thinspace\thinspace\thinspace b, |b\thinspace\thinspace\thinspace a, a|b, b|a, \}$, $| Q^1_2|$=$J^1_2$(by definition of $Q_n^1$). We say $Q_1^1=\{|a, a|\}$. In proving the corollary we generalise the method used in \cite{barred:2013} in proving theorem 2.1 from a single label of a bar to a double label of a bar, in the following way.  We consider a map $f:\{0,1\}_2\times\{0,1\}_1\times Q_1^1\rightarrow Q^1_2$. We construct the set$\{0,1\}_2\times\{0,1\}_1\times Q_1^1$ in the following way, \\ 1.  $\{0,1\}_1$ labels one of the elements of $Q_1^1$ with a binary label 0 or 1 at the bottom. Hence from the two element of $Q^1_1$ there will result four labeled barred preferential arrangements.  The set $\{0,1\}_1\times Q^1_1=\{\underset{0}{|}a, \underset{1}{|}a, a\underset{0}{|}, a\underset{1}{|}\}$.
\\2. Then $\{0,1\}_2$ labels at the top, the four elements of $\{0,1\}_1\times Q^1_1$ to form eight elements of  $\{0,1\}_2\times\{0,1\}_1\times Q_1^1$ in the following way,\\ 
 $\{0,1\}_2\times\{0,1\}_1\times Q^1_1=\{ \overset{0}{\underset{0}{|}} a,\thinspace
 a\overset{0}{\underset{0}{|}},\thinspace \overset{0}{\underset{1}{|}}a,\thinspace  a\overset{0}{\underset{1}{|}}, \thinspace  
 \overset{1}{\underset{0}{|}} a,\thinspace
 a\overset{1}{\underset{0}{|}},\thinspace \overset{1}{\underset{1}{|}}a,\thinspace  a\overset{1}{\underset{1}{|}}
\} $
We now define a mapping from $\{0,1\}_2\times\{0,1\}_1\times Q^1_1$ to $Q^1_2$ in the following way,\\ 
\RNum{1}. If the indexes on the bar are both 0, then for a barred preferential arrangement, by introducing a second element $b$, place $a$ and $b$ as a single block on the side of the bar where $a$ was(in all cases remove the indices on the bars). An example of such an arrangement is $|ab$.
\\\RNum{2}. If the top index on the bar is 0 and the bottom index on the bar is 1, then remove the indexing on the bar and place $b$ on opposite side of the bar from $a$ in obtaining an element of $Q^1_2$. An example of such an arrangement is $a|b$.
\\\RNum{3}. If the top index on the bar is 1 and the bottom index is 0, then place the two elements on the same side of the bar but placing the element $b$ on its on block immediately before the bar then following by the block of the element $a$. An example of such an arrangement is,  $a\thinspace\thinspace\thinspace b|$. 
\\\RNum{4}. If both index are equal to 1, then place both elements on separate block on the same side of the ba, where now as opposed to 3 above the immediate block before the bar will be that of $a$ then followed by the block of $b$. an example of such an arrangement is $|a\thinspace\thinspace\thinspace b$.
\\From the way we have constructed the function $f$,
\\two elements of $\{0,1\}_2\times\{0,1\}_1\times Q_1^1$ satisfy each of the four conditions above and each element of $\{0,1\}_2\times\{0,1\}_1\times Q_1^1$ is mapped to one element of $Q^1_2$. Hence the mapping is 1-1.
\\ Clear by marking the position of the elements $a$ and $b$ in an element $\textgoth{X}\in Q^1_2$, we can identify an element of  $\{0,1\}_2\times\{0,1\}_1\times Q_1^1$ which is mapped to $\textgoth{X}$. Hence the map in onto. So $|\{0,1\}_2\times\{0,1\}_1\times Q_1^1|=|Q^1_2|\Rightarrow 2\times2\times J^1_1=J^1_2$ as required. 
\begin{corollary}for $m=1$ and $n=1$ the identity holds,
\begin{center}$J^1_1=J^1_0+J^0_1$\end{center}
\end{corollary}
We prove the corollary by considering a mapping $f:Q^1_1\rightarrow Q^1_0\cup Q^0_1$, where $Q^m_n$ is the set of all barred preferential arrangements of an $n\geq1$ element set having $m\geq0$ bars. The set $Q^1_1=\{|a, a|\}$, also $Q^1_0\cup Q^0_1=\{|, a\}$. We define a mapping from $Q^1_1$ in the following way,
\\\RNum{1}. If a is before the bar on $\textgoth{X}\in Q^1_1$ then map $\textgoth{X}$ to $a$ in $Q^1_0\cup Q^0_1$ (this operation behaves as if you remove the bar from $\textgoth{X}$, when mapping $\textgoth{X}$ to an element of $Q^1_0\cup Q^0_1$).
\\\RNum{2}. If $a$ is after the bar on $\textgoth{X}\in Q^1_1$ then map $\textgoth{X}$ to $|$ in $Q^1_0\cup Q^0_1$ (this operation behaves as if $a$ is removed from $\textgoth{X}$ when mapping $\textgoth{X}$ to an element of   $Q^1_0\cup Q^0_1$). Clearly this mapping is a one to one and onto mapping between the two sets of two elements each. So $|Q^1_1|=|Q^1_0\cup Q^0_1|$. Hence $J^1_1=J^0_1+J^1_0$ ($|Q^1_0\cup Q^0_1|=J^0_1+J^1_0$ this is so, since the two sets are disjoint). Thus the corollary. 
\begin{corollary}
for $m=1$ and $n=2$ the identity holds,
\begin{center}
$J^1_2=3\times J^1_1+2\times J^1_0$
\end{center}
\end{corollary}

We prove the corollary by establishing a bijection between a set $Q^1_2$ whose cardinality is $J^1_2$ and a set $\{0,1,2\}\times Q^1_1\cup\{0,1\}\times Q^1_0$, whose cardinality is $3\times J^1_1+2\times J^1_0$. Where $Q^m_n$ denotes the set of all barred preferential arrangements of an $n\geq1$ \\element set having $m\geq0$ bars. We say, $Q_1^1=\{|a, a|\}$. We then construct the set $\{0,1,2\}\times Q^1_1$ in the following way, $\{0,1,2\}\times Q^1_1=\{a\underset{0}{|},\thinspace\thinspace \underset{0}{|}a, \thinspace\thinspace a\underset{1}{|},\thinspace\thinspace \underset{1}{|}a,\thinspace\thinspace a\underset{2}{|},\thinspace\thinspace \underset{2}{|}a\}$,So the bar of each element of $Q^1_1$ is labeled by either 0, 1 or 2, in forming the set $\{0,1,2\}\times Q^1_1$. We construct the set $\{0,1\}\times Q^1_0$ in a similar way to obtain
$\{0,1\}\times Q^1_1=\{\underset{0}{|},\thinspace\thinspace \underset{1}{|}\}$. So we have, $\{0,1,2\}\times Q^1_1\cup\{0,1\}\times Q^1_0=\{a\underset{0}{|},\thinspace\thinspace \underset{0}{|}a, \thinspace\thinspace a\underset{1}{|},\thinspace\thinspace \underset{1}{|}a,\thinspace\thinspace a\underset{2}{|},\thinspace\thinspace \underset{2}{|}a,\thinspace\thinspace \underset{0}{|},\thinspace\thinspace \underset{1}{|} \}$   
. Also we have, 
$Q^1_2=\{a|b, b|a,
\thinspace\thinspace  |a\thinspace\thinspace b, \thinspace\thinspace
 b\thinspace\thinspace a|,\thinspace\thinspace a \thinspace\thinspace b|
 ,\thinspace\thinspace |b\thinspace\thinspace\thinspace a, \thinspace\thinspace|ab,\thinspace\thinspace ab|\}
 $.
 We define a mapping $f:\{0,1,2\}\times Q^1_1\cup\{0,1\}\times Q^1_0\rightarrow Q^1_2$ in the following way,
\\ \RNum{1}. If $\textgoth{X}\in \{0,1,2\}\times Q^1_1$ and the indexing on the bar being 0. Then place an element $b$ on opposite side of $a$. An example of such an arrangement is $a|b$. 
\\\RNum{2}. If $\textgoth{X}\in \{0,1,2\}\times Q^1_1$ and the indexing on the bar being 1, then place an element $b$ on it's own block on the same side as $a$ but the block of $a$ should be the one which is immediately before the bar. An example of such an arrangement is $|a\thinspace\thinspace b$.
\\\RNum{3}. If $\textgoth{X}\in \{0,1,2\}\times Q^1_1$ and the indexing on the bar being 2. Then place an element $b$ on its own block immediately before the bar, then followed by the block of $a$. An example of such an arrangement is $|b\thinspace\thinspace\thinspace a$.
\\\RNum{4}.If $\textgoth{X}\in \{0,1\}\times Q^1_0$ and the index on the bar is 0 then introduce write two elements $a$ and $b$ as a single block to the right of the bar. Hence that element will be mapped to the element $|ab$ in $Q^1_2$. 
\\\RNum{5}. If $\textgoth{X}\in \{0,1\}\times Q^1_0$ and the index on the bar is 1 then introduce write two elements $a$ and $b$ as a single block to the left of the bar. Hence that element will be mapped to the element $ab|$ in $Q^1_2$.      
\\From the way we have defined $f$ above, the $1^{st}$ and the $2^{nd}$ elements of $Q^1_2$, satisfy property \RNum{1}. The $3^{rd}$ and the $4^{th}$ element of $Q^1_2$ satisfy property \RNum{2}. The $5^{th}$ and the $6^{th}$ elements of $Q^1_2$ satisfy property \RNum{3}. The $7^{th}$ element of $Q^1_2$ satisfy property \RNum{4} above. The $8^{th}$ element of $Q^1_2$ satisfy property \RNum{5} above. From the way we have defined $f$ above, each element of $\{0,1,2\}\times Q^1_1\cup\{0,1\}\times Q^1_0$ is mapped to a unique element of $Q^1_2$. Also using $f$ above, we can move from elements of $Q^1_2$ to elements of $\{0,1,2\}\times Q^1_1\cup\{0,1\}\times Q^1_0$ with out ambiguity. Hence the mapping is bijective. Thus $|\{0,1,2\}\times Q^1_1\cup\{0,1\}\times Q^1_0|=|Q^1_2|$ i.e $3\times J^1_1+2\times J^1_0=J^1_2$ as required(note the sets $\{0,1,2\}\times Q^1_1$ and the set $\{0,1\}\times Q^1_0$ are disjoint hence order of their union is the sum of the orders).  

\begin{theorem} for $m,n\geq1$,
 \begin{center}$J^k_n=J^{k-1}_{n}+\sum\limits_{s=0}^{n-1}\binom{n}{s}J^{k-1}_{s}J^{0}_{n-s}$\end{center}
 \end{theorem}
 We recall  $J^k_n$ is the cardinality of a set $Q^k_n$ which is the collection of all possible barred preferential arrangements of an $n$-element having $k$ bars. Each barred preferential arrangement in $Q^k_n$ has $k+1$ section since there are $k$ bars. We prove the theorem by partitioning $Q^k_n$ into disjoint subsets $\mu_1$ and $\mu_2$. Where $\mu_1$ is the collection of all those elements from $Q^k_n$ whose first section is empty and $\mu_2$ is the collection of all those elements from $Q^k_n$ whose first section is non-empty. In finding the number of element in $\mu_1$ we argue as follows: the first section of $\textgoth{O}\in Q^k_n$ being empty means the $n$ elements are distributed among the $k$ other sections. A distribution of $n$ elements among $k$ sections is a distribution of the elements among $k-1$ bars. As a result $|\mu_1|=J^{k-1}_{n}$.  
 
 For the case the first section of each element of $Q^k_n$ being required to be non-empty we argue as follows in finding the number of elements of $\mu_2$: There can be a minimum of 0 elements not in the first section and there can be a maximum of $n-1$ elements which are not in the first section on all barred  preferential arrangements  which are to be in $\mu_2$. Lets assume there are $s$ elements which are not in the first section of each  barred preferential arrangement. There are $\binom{n}{s}$ ways of selecting the $s$ elements. There are $J^{k-1}_{s}$ ways of preferentially arranging the $s$ elements among the other $k$ sections. There $n-s$ remaining elements can be preferentially arranged on the first section in $J^0_{n-s}$ ways. Taking the product and summing over $s$ we obtain $\mu_2=\sum\limits_{s=0}^{n-1}\binom{n}{s}J^{k-1}_{s}J^0_{n-s}$. When we combine $\mu_1$ and $\mu_2$ we obtain the result of the theorem.

 %In \cite{nel:1991} the authors proposed the following identity. 
 % \begin{equation}\label{equation:1}2p_{n}(t)-p_{n}(t+1)=t^n\qquad where\thinspace\thinspace\thinspace t,n\geq0\end{equation} 
\section{Identities of restricted barred preferential arrangements}
 We denote by $H_n^m$ the set of all barred preferential arrangements of an $n$ element set having $m$ bars, in which $m$  fixed  sections (out of $m+1$ sections) are  allowed to have a maximum of only one  block and the other section can have one or more blocks (we refer to this section as the free  section). The elements of the free section can have any possible number of blocks of a preferential arrangement of given $r\geq0$ elements. So $H_n^m$ represents the set of all barred preferential arrangements of an $n$ element set having $m$ bars in which the $m$ fixed (first) sections have a maximum of one block but one fixed section can have more than one block. The free section can be any of the section as long its fixed.  We denote by $|H_n^m|=I_n^m$.

 We prove the identity in the following way,
 \begin{lemma}\label{lemma:2}\thinspace
 for\thinspace\thinspace\thinspace $m=0$ and $n\geq0$ we have,
 \begin{center}
$ 2I_n^0 = I_n^1$
 \end{center}
 \end{lemma}

%TO WORK FROM HERE LATER

 We consider the set $H_n^0$ which is the set of all barred preferential arrangements of an $n$ element set with no bars. Hence all elements of $H_n^0$  have only one section since there are no bars. We also consider the set $H_n^1$, which is the set of all barred  preferential arrangement of an $n$ element set having one bar, in which the first section has a at most one block and the second section is just a preferential arrangement of given elements not necessarily into one block.\\

Adding an extra bar $\overset{*}{{|}}$ to the far right of each element of the set $H_n^0$ does not affect counting. We add an extra bar $\overset{*}{{|}}$ to the far right of each element of the set $H_n^0$ to form the set $*\P ^0_n$, with the extra bar on the far right of each element.In proving the lemma, use a similar method to the one used in \cite{barred:2013} in proving theorem~2.1
. We have $|*\P ^0_n|= I_n^0 = |H_n^0|$. We define the set  $R^*_n(0)=\{0,1\}\times *\P ^0_n$ as containing the same elements as  $*\P ^0_n$ but with an indexing on the bar $\overset{*}{{|}}$ which is either 0 or 1, hence $R^*_n(0)$ has twice the number of elements as $*\P ^0_n$ (half of them having the index 0 and the other half the index 1). 
We now use the elements of $R^*_n(0)$ to reconstruct the set $H_n^1$. We do it as follows:\\

%We are here 21 10 2014

\RNum{1}. If the indexing on the bar $\overset{*}{{|}}$ on an element $\textgoth{X}\in R^*_n(0)$ is 0 then such an element will be interpreted in $H^{1}_n$, as an element of  $H^{1}_n$ whose ${2}^{nd}$ section is empty. We collect all such elements to form the set $W$. The set $W$ has $I_n^0$ elements(half of the elements in $R^*_n(0)$ have the index 0)
% If the indexing on the bar $\overset{*}{{|}}$ on an element $x\in R^*_n(0)$ is 0 then such an element will be interpreted in $H_n^1$, those elements in  $H_n^1$ whose ${2}^{nd}$ section is empty. 
. 
\\\RNum{2} If the indexing on the bar $\overset{*}{{|}}$ on an element $\textgoth{X}\in R^*_n(0)$ is 1 then we shift the last block of section 1 to be the only block to the right of the bar $\overset{*}{{|}}$ to form the set $K$. The set $K$ has $I_n^0$ elements(half of the elements in $R^*_n(0)$ have the indexing 1)\\
Clearly the sets $K$ and $W$ are disjoint. So we have $|K\cup W|=|K|+|W| \Rightarrow |K\cup W|=I_n^0+I_n^0$. Now what we are having is two sets $K\cup W$ and $H^{1}_n$ have the same definition(one fixed section can have a maximum of one block and the other section is a free section) and the same number of elements  in the underlying sets, So  they must be of the same size. Hence $|K\cup W|=I_n^1+I_n^1= I_n^1 = |H^{1}_n|$. Thus the lemma.\\  

For the value $m=1$, we have the following lemma,
\begin{lemma}\label{lemma:3}
\thinspace
 for\thinspace\thinspace\thinspace $m=1$ and $n\geq0$ we have,
 \begin{center}
$ 2I_n^1-1=I_n^2$
 \end{center}
\end{lemma} 
 We consider a set $H^1_n$, which is the set of all barred preferential arrangements of an $n$ element set having $1$ bar, in which the first section has a maximum of one block each and the ${2}^{nd}$ section to be just a preferential arrangement of  given elements not necessarily into one block (this is the only section with this property). So the ${2}^{nd}$ section of each element of $H^1_n$ is a free section. We also consider the set $H^{2}_n$, which  is the set of all barred preferential arrangements of an $n$-element set, having $2$ bars in which  $2$ fixed sections have a maximum of one block each and one section can have more than one block (the free section).

 We now want to reconstruct $H^{2}_n$ using elements of $H^1_n$. Adding an extra bar $\overset{*}{{|}}$ to the far right of each element of $H^1_n$ does not affect counting. We do that to form the set $*\P ^1_n$. Now we have $|*\P ^1_n|=I_n^1=|H^1_n|$. The bar $\overset{*}{{|}}$ is to the right of the free sections of elements of $H^1_n$. On each element $\textgoth{X}\in*\P_{n}^1$ to the left of the bar $\overset{*}{{|}}$ is the free section  of $\textgoth{X}$ and the section to the right of $\overset{*}{{|}}$ is empty.
 \\     
  We define the set  $R^*_n(1)=\{0,1\}\times *\P ^1_n$ as containing the same elements as  $*\P ^1_n$ but with an indexing on the bar $\overset{*}{{|}}$ which is either 0 or 1, hence $R^*_n(1)$ has twice the number of elements as $*\P ^1_n$(of which half of them having the indexing 0 and the other half having the index 1).    
We now use elements of $R^*_n(1)$ to reconstruct the set $H^{2}_n$. We construct as follows,
\\\RNum{1}.If the indexing on the bar $\overset{*}{{|}}$ on an element $\textgoth{X}\in R^*_n(1)$ is 0 then such an element will be interpreted in $H^{2}_n$, as an element of  $H^{2}_n$ whose ${3}^{rd}$ section is empty.
 %If the indexing on the bar $\overset{*}{{|}}$ on an element $x\in R^*_n(7)$ is 0 then such elements will be interpreted in $P^{8}_n$, as those elements in  $P^{8}_n$ whose ${9}^{th}$ section is empty.
  We collect all such elements to form the set $W$. The set $W$ has $I_n^1$ elements(half of the elements in $R^*_n(1)$ have the indexing 0)
\\\RNum{2} If the indexing on the bar $\overset{*}{{|}}$ on an element $\textgoth{X}\in R^*_n(1)$ is 1 then we shift the last block of the ${3}^{rd}$ section of $\textgoth{X}$ to be the only block to the right of $\overset{*}{{|}}$(i.e the block closest to the bar $\overset{*}{{|}}$ on $\textgoth{X}$) to form the set $K$. There are $I_n^1$ elements having index 1 in $R^*_n(1)$(half of the elements in $R^*_n(1)$ have the indexing 1). On the construct of the set $K$ some elements of $K$ also appear as elements of the set $W$. Those redundant elements(elements having identical interpretation in $H^{2}_n$ ) occur when there free section on elements of $R^*_n(1)$ is empty. In that case there is no block to put to the right of the bar $\overset{*}{{|}}$ when constructing the set $K$. In that case we get common elements between $K$ and $W$. Common elements between $K$ and $W$ occur when the $n$ elements are distributed on the first $m$ sections of $H^1_n$ of which there is only one way of doing that. So there is one  element in the intersection of $K$ and $W$(when elements of $K$ and $W$ are interpreted as elements of ).

The elements of  both sets $K$ and $W$ have $3$ sections, of which for those with a fixed $\textgoth{X}\in K\cup W$,  $2$ sections of $\textgoth{X}$ have a maximum of one block and one section is a free section. We have, $|K\cup W|=|K|+|W|-|K\cap W|$$\Rightarrow$ $|K\cup W|=I_n^1+I_n^1-1$. Now what we are having is two sets $K\cup W$ and $H^{2}_n$ have the same definition and the same number of elements in the underlying sets, hence must be of the same size. That is, $|K\cup W|=I_n^1+I_n^1-1 = I_n^2 = |H^{2}_n|$. Thus the lemma.\\

The case $m=2$ differs from the above two cases in the sense that a non-trivial constant term arises in the form of $2^n$. 
Therefore we state and prove for the value $m=2$, the following lemma,
\begin{lemma}\label{lemma:4}
\thinspace
 for\thinspace\thinspace\thinspace $m=2$ and $n\geq0$ we have,
 \begin{center}
$ 2I_n^2-2^n=I_n^3$
 \end{center}
\end{lemma} 
 We consider a set $H^2_n$, which is the set of all barred preferential arrangements of an $n$ element set having $2$ bars, in which the first $2$ sections have a maximum of one block each and the ${3}^{rd}$ section to be just a preferential arrangement of  given elements not necessarily into one block( this is the only section with this property). So the ${3}^{rd}$ section of each element of $H^2_n$ is the free section. We also consider the set $H^{3}_n$, which  is the set of all barred preferential arrangements of an $n$ element set, having $3$ bars in which  $3$ fixed sections have a maximum of one block each and one section can have more than one block( the free section).

 We now want to reconstruct $H^{3}_n$ using elements of $H^2_n$. Adding an extra bar $\overset{*}{{|}}$ to the far right of each element of $H^2_n$ does not affect counting. We do that to form the set $*\P ^2_n$. Now we have $|*\P ^2_n|=I_n^2=|H^2_n|$.  Where the bar $\overset{*}{{|}}$ is to the right of the free sections of elements of $H^2_n$. On each element $\textgoth{X}\in*\P_{n}^2$ to the left of the bar $\overset{*}{{|}}$ is the free section  of $\textgoth{X}$ and the section to the right of $\overset{*}{{|}}$ is empty.
 \\     
  We define the set  $R^*_n(2)=\{0,1\}\times *\P ^2_n$ as containing the same elements as  $*\P ^2_n$ but with an indexing on the bar $\overset{*}{{|}}$ which is either 0 or 1, hence $R^*_n(2)$ has twice the number of elements as $*\P ^2_n$(of which half of them having the indexing 0 and the other half having the index 1).    
We now use elements of $R^*_n(2)$ to reconstruct the set $H^{3}_n$. We construct as follows,
\\\RNum{1}.If the indexing on the bar $\overset{*}{{|}}$ on an element $\textgoth{X}\in R^*_n(2)$ is 0 then such an element will be interpreted in $H^{8}_n$, as an element of  $H^{3}_n$ whose ${4}^{th}$ section is empty.
 %If the indexing on the bar $\overset{*}{{|}}$ on an element $x\in R^*_n(7)$ is 0 then such elements will be interpreted in $P^{8}_n$, as those elements in  $P^{8}_n$ whose ${9}^{th}$ section is empty.
  We collect all such elements to form the set $W$. The set $W$ has $I_n^2$ elements(half of the elements in $R^*_n(2)$ have the indexing 0)
\\\RNum{2} If the indexing on the bar $\overset{*}{{|}}$ on an element $\textgoth{X}\in R^*_n(2)$ is 1 then we shift the last block of the ${3}^{rd}$ section of $\textgoth{X}$ to be the only block to the right of $\overset{*}{{|}}$(i.e the block closest to the bar $\overset{*}{{|}}$ on $\textgoth{X}$) to form the set $K$. There are $I_n^2$ elements having index 1 in $R^*_n(2)$(half of the elements in $R^*_n(2)$ have the indexing 1). On the construct of the set $K$ some elements of $K$ also appear as elements of the set $W$. Those redundant elements occur when there free section on elements of $R^*_n(2)$ is empty. In that case there is no block to put to the right of the bar $\overset{*}{{|}}$ when constructing the set $K$. In that case we get common elements between $K$ and $W$. Common elements between $K$ and $W$ occur when the $n$ elements are distributed on the first $m$ sections of $H^2_n$ of which there are $2^n$ ways of doing that. So the number of elements in the intersection of $K$ and $W$ is $2^n$.

The elements both sets $K$ and $W$ have $4$ sections. Of which for a fixed $\textgoth{X}\in K\cup W$,  $3$ sections of $\textgoth{X}$ have a maximum of one block and one section is a free section. We have, $|K\cup W|=|K|+|W|-|K\cap W|$$\Rightarrow$ $|K\cup W|=I_n^2+I_n^2-2^n$. Now what we are having is two sets $K\cup W$ and $H^{3}_n$ have the same definition and the same number of elements  in the underlying sets, So  they must be of the same size. Hence $|K\cup W|=I_n^2+I_n^2-2^n = I_n^3 = |H^{3}_n|$. Thus the lemma.

  We generalise lemmas \ref{lemma:2},\ref{lemma:3} and \ref{lemma:4}  into the following theorem,
\begin{theorem}\label{theorem:2}
\thinspace
 for\thinspace\thinspace\thinspace $m,n\geq0$, we have
 \begin{center}
$ 2I_n^m-m^n = I_n^{m+1}$
 \end{center}
\end{theorem}   
%The proof of theorem \ref{theorem:2} follows from that of lemma \ref{lemma:2} and lemma \ref{lemma:3}.  
 
 We consider a set $H^m_n$, which is the set of all barred preferential arrangements of an $n$ element set having $m$ bars, in which the first $m$ sections have a maximum of one block each and the ${(m+1)}^{th}$ section to be just a preferential arrangement of  given elements not necessarily into one block( this is the only section with this property). So the ${(m+1)}^{th}$ section of each element of $H^m_n$ is the free section. We also consider the set $H^{m+1}_n$, which  is the set of all barred preferential arrangements of an $n$ element set, having $m+1$ bars in which  $m+1$ fixed sections have a maximum of one block each and one section can have more than one block(the free section).

 We now want to reconstruct $H^{m+1}_n$ using elements of $H^m_n$. Adding an extra bar $\overset{*}{{|}}$ to the far right of each element of $H^m_n$ does not affect counting. We do that to form the set $*\P ^m_n$. Now we have $|*\P ^m_n|=I_n^m=|H^m_n|$.  Where the bar $\overset{*}{{|}}$ is to the right of the free sections of elements of $H^m_n$. On each element $\textgoth{X}\in*\P_{n}^m$ to the left of the bar $\overset{*}{{|}}$ is the free section  of $\textgoth{X}$ and the section to the right of $\overset{*}{{|}}$ is empty.
 \\     
  We define the set  $R^*_n(m)=\{0,1\}\times *\P ^m_n$ as containing the same elements as  $*\P ^m_n$ but with an indexing on the bar $\overset{*}{{|}}$ which is either 0 or 1, hence $R^*_n(m)$ has twice the number of elements as $*\P ^m_n$(of which half of them having the indexing 0 and the other half having the index 1).    
We now use elements of $R^*_n(m)$ to reconstruct the set $H^{m+1}_n$. We construct as follows,
\\\RNum{1}.If the indexing on the bar $\overset{*}{{|}}$ on an element $\textgoth{X}\in R^*_n(m)$ is 0 then such an element will be interpreted in $H^{m+1}_n$, as an element of  $H^{m+1}_n$ whose ${(m+2)}^{th}$ section is empty.
% If the indexing on the bar $\overset{*}{{|}}$ on an element $x\in R^*_n(t)$ is 0 then such elements will be interpreted in $P^{t+1}_n$, as those elements in  $P^{t+1}_n$ whose ${(t+2)}^{th}$ section is empty.
 We collect all such elements to form the set $W$. The set $W$ has $I_n^m$ elements(half of the elements in $R^*_n(m)$ have the indexing 0)
\\\RNum{2} If the indexing on the bar $\overset{*}{{|}}$ on an element $\textgoth{X}\in R^*_n(m)$ is 1 then we shift the last block of the ${(m+1)}^{th}$ section of $\textgoth{X}$ to be the only block to the right of $\overset{*}{{|}}$(i.e the block closest to the bar $\overset{*}{{|}}$ on $\textgoth{X}$) to form the set $K$. There are $I_n^{m}$ elements having index 1 in $R^*_n(m)$(half of the elements in $R^*_n(m)$ have the indexing 1). On the construct of the set $K$ some elements of $K$ also appear as elements of the set $W$. Those redundant elements occur when there free section on elements of $R^*_n(m)$ is empty. In that case there is no block to put to the right of the bar $\overset{*}{{|}}$ when constructing the set $K$. In that case we get common elements between $K$ and $W$. Common elements between $K$ and $W$ occur when the $n$ elements are distributed on the first $m$ sections of $H^m_n$ of which there are $m^n$ ways of doing that. So the number of elements in the intersection of $K$ and $W$ is $m^n$.

The elements both sets $K$ and $W$ have $m+2$ sections. Of which for a fixed $\textgoth{X}\in K\cup W$,  $m+1$ sections of $\textgoth{X}$ have a maximum of one block and one section is a free section. We have, $|K\cup W|=|K|+|W|-|K\cap W|$$\Rightarrow$ $|K\cup W|=I_n^m+I_n^m-m^n$. Now what we are having is two sets $K\cup W$ and $H^{m+1}_n$ have the same definition and the same number of elements  in the underlying sets, So  they must be of the same size. Hence $|K\cup W|=I_n^m+I_n^m-m^n = I_n^{m+1}=|H^{m+1}_n|$. Thus the theorem.
\\
\section{Acknowledgments}
Both authors acknowledge the support from Rhodes University. The first author would like also to acknowledge financial support from the DAAD-NRF scholarship of South Africa, the Levenstein Bursary of Rhodes University and the NRF-Innovation doctoral scholarship of South Africa.


\begin{thebibliography}{10}
\bibitem{cayley:1859}
A. Cayley,\emph{ LVIII. On the analytical forms called trees.–Part II}, The London, Edinburgh, and Dublin Philosophical Magazine and Journal of Science 18,121 (1859): 374-378.
\bibitem{com:74} L. Comtet, \emph{Advanced Combinatorics}, Reidel, 1974, p. 294.
\bibitem{gross:1962} 
O. A. Gross, \emph{Preferential arrangements}, Amer. Math. Monthly, 69 (1962) 4-8.
\bibitem{mendelson:1982}
Mendelson Elliott,
\emph{Races with Ties}, Mathematics Magazine, Vol.55, No.3(May,1982), 170-175.
\bibitem{murali:Nelsongenerating function} V. Murali, \emph{Ordered Partitions and Finite Fuzzy Sets}, Far East J. Math.Sci.(FJMS), 21(2) (2006), 121-132.
\bibitem{oeis} https://oeis.org/,  The OEIS Foundation, 1964.
%\bibitem{nel:91} R.B.Nelsen and H.Schmidt Jr., Chains in power sets,  Math.Mag., 64 (1991) 23-31.
\bibitem{Pippenger:2010}
 Pippenger  Nicholas,
\emph{The Hypercube of resistors, Asymptotic Expansions and Preferential Arrangements},
The American mathematical monthly,83(2010) pp.331-346.
\bibitem{barred:2013}
C. Ahlbach, J. Usatine, and N. Pippenger. \emph{Barred Preferential Arrangement}, Electronic Journal of Combinatorics, 20:2 (2013), 1-18
 
\end{thebibliography}
\end{document}